\documentclass[9 pt]{amsart} 

\usepackage{
amsfonts,
latexsym,
amssymb,
enumerate,
verbatim,
mathrsfs,
graphicx,
centernot,
color,
bbding}

\newcommand{\myplus}{ \scriptscriptstyle \textup{\FiveStarOutlineHeavy}}

\newtheorem{theorem}{Theorem}[section]

\newtheorem{proposition}[theorem]{Proposition} 
 
\newtheorem{corollary}[theorem]{Corollary}

\newtheorem*{claim*}{Claim}

\newtheorem*{theorem*}{Theorem}
\newtheorem*{proposition*}{Proposition}
\newtheorem*{corollary*}{Corollary}
\newtheorem*{lemma*}{Lemma}
\newtheorem*{scholion*}{Scholion}

\theoremstyle{definition}
\newtheorem{definition}[theorem]{Definition}

\theoremstyle{remark}
\newtheorem{remark}[theorem]{Remark}

\newtheorem*{remark*}{Remark}
\newtheorem*{remarks*}{Remarks}
 
\newtheorem{example}[theorem]{Example}

\newtheorem*{observation*}{Observation}

 \allowdisplaybreaks[1]

\numberwithin{equation}{section}

\begin{document}

\title{Semilattices with a congruence}

\author{Paolo Lipparini} 
\address{Dipartimento di Matematica\\Viale della  Ricerca
Relazionale\\Universit\`a di Roma ``Tor Vergata'' 
\\I-00133 ROME ITALY (Currently retired)
\\ORCiD 0000-0003-3747-6611}

\email{lipparin@axp.mat.uniroma2.it}



\thanks{Work performed under the auspices of G.N.S.A.G.A. 
The author acknowledges the MIUR Department Project awarded to the
Department of Mathematics, University of Rome Tor Vergata, CUP
E83C18000100006.}

\begin{abstract}
A specialization semilattice is a semilattice together 
with a coarser preorder satisfying a compatibility condition.
We show that the category of
specialization semilattices is isomorphic
to the category of semilattices with a congruence,
hence equivalent to the   
 category of  semilattice epimorphisms. 

Guided by the above  example,
we recall an ``internal'' characterization of surjective
homomorphisms between general relational systems.
 \end{abstract}

\maketitle

\section{Introduction} \label{intro} 

A \emph{specialization semilattice}  is a join semilattice
endowed with a further preorder $ \sqsubseteq $ which is coarser than
the order induced by the semilattice operation and such that 
\begin{align}
\label{s3}  
 &a \sqsubseteq b  \ \text{ and } \   a_1 \sqsubseteq b
\ \text{ imply } \ 
 a \vee a_1 \sqsubseteq b, 
\end{align}
for all the elements $a,a_1,b$
of the semilattice. 

If $X$ is a topological space, then
$( \mathcal P(X),   {\cup}, {\sqsubseteq} )$
is a specialization semilattice, 
where $ \sqsubseteq $ is
defined by
$ a \sqsubseteq b$
if  $a \subseteq  K b$,
for $a$ and $b$ subsets of $X$,  
and $ K $ is closure.
In \cite[Theorem 5.7]{mtt} we showed that every
specialization semilattice can be embedded  
into the specialization semilattice associated
to some topological space,  as described above.

The study of specialization semilattices
is motivated by the fact that a function between topological spaces
is continuous if and only if the corresponding image
function is a homomorphism 
between the associated specialization semilattices \cite[Proposition 2.4]{mtt}.
Alternative ``algebraizations'' of topology do not satisfy
such an exact functorial condition. See the introduction
of \cite{mtt} for a discussion. 

More or less in disguise and
in different terminology,
specialization semilattices appear  in domain theory, in 
algebraic geometry, in lattice theory,
in set theory (for example if we take inclusion modulo finite
as a ``specialization''), in algebraic logic, in the study of tolerance spaces,
 in measure theory and as \emph{Complete Implicational Systems},
motivated by  relational data bases,
data analysis and  artificial intelligence.
See \cite[Section 4]{mtt} for details and references. 

Besides the  topological construction hinted above, 
there is another way of representing
specialization semilattices.
If $(S, \vee_S)$,  $(T, \vee_T)$ are semilattices,
$\varphi$  is a homomorphism from $\mathbf S$ to $\mathbf T$
and in $S$ we set 
\begin{equation}\label{sqf}    
\text{$a \sqsubseteq _ \varphi b$
if $\varphi(a) \leq _T \varphi (b)$,} 
  \end{equation} 
then  $(S, \vee_S, \sqsubseteq _ \varphi )$  
is a specialization semilattice.
In the formula \eqref{sqf}  $\leq_T$ is the partial order on 
$T$ induced by the semilattice operation $\vee_T$. 
 
In Corollary \ref{rep} below 
we will show that every specialization semilattice can be represented
in the above way.   
Thus specialization semilattices are at the same time
``substructures'' of topological spaces and ``kernels''
of semilattice homomorphisms.
Together with the mentioned abundance of 
examples occurring in varied mathematical contexts,
this double kind of representation confirms
the naturalness of the notion.

We will also show that, given a semilattice $\mathbf S$,
the set of preorders making $\mathbf S$ a
specialization semilattice is in a bijective correspondence
with the set of congruences on $\mathbf S$;
actually, the correspondence is a complete lattice
isomorphism. See Theorem \ref{latt}.
Put in another way, the set of those 
preorders which make $\mathbf S$ a
specialization semilattice is in  a bijective
correspondence with the 
class of surjective homomorphisms
with domain $\mathbf S$, considered up
to isomorphisms keeping $\mathbf S$ fixed. 

As another reformulation,  the category of
specialization semilattices, with the standard notion of
morphism, is isomorphic to the category of
semilattices with a congruence, where in the latter
category morphisms
are semilattice morphisms which respect the
congruence. See Corollary \ref{corcat}. It is then standard to see that 
the category of specialization semilattices
is equivalent to the category  $\mathcal S\mathfrak{Epi}$ 
of semilattice epimorphisms. Objects of 
$\mathcal S\mathfrak{Epi}$ are semilattice epimorphisms
and morphisms of $\mathcal S\mathfrak{Epi}$ are commuting squares
of semilattice morphisms,
where the additional  morphisms in the square
are not necessarily assumed to be epi.

Let us now consider  a weaker notion. A \emph{specialization poset} 
is a partially ordered set  (henceforth, poset, for short) together 
with a coarser preorder. Again, examples of specialization
posets abound in the mathematical literature \cite[Section 4]{mtt}. 
The connection with topological spaces
goes the same way as for specialization semilattices.
If $X$ is a topological space, then
$( \mathcal P(X),   { \subseteq }, {\sqsubseteq} )$ 
is a specialization poset; every
specialization poset can be represented as a substructure
of such a topological specialization poset. See \cite[Proposition 5.10]{mtt}. 
As for ``quotient'' representations, if $(P, \leq _P)$, 
 $(Q, \leq _Q)$ are posets and $\varphi$  is an order preserving
map, then we can define $ \sqsubseteq _ \varphi $ as in 
\eqref{sqf}. Then  $(P, \leq _P, \sqsubseteq _ \varphi)$
turns out to be a specialization poset. Every
specialization poset can be represented this way.
See Example 3.1.
As above,  the category of specialization posets
is equivalent to the category $\mathcal{POS}\mathfrak{Surj}$ of  
 surjective order preserving
maps between posets, where, again, morphisms of 
$\mathcal{POS}\mathfrak{Surj}$
are commuting squares and the additional poset morphisms
are not assumed to be surjective. 
However, it is standard to see that surjective order preserving
maps between posets cannot be
characterized ``internally'' by an equivalence relation,
see Remark \ref{diff}.  

Hence an ``internal''  characterization of 
surjective order preserving
maps between posets is better provided by a
binary relation (a coarser 
preorder) which is not necessarily
an equivalence relation. Put in another way,
specializations in posets can be used exactly
as congruences in algebraic structures,
where we have only operations and no relation,
apart from equality.

In Section \ref{gen} we recall a similar characterization
for models for an arbitrary language $\mathscr L$,
allowed to contain both relation and function symbols.
In this sense, the ``kernel'' of some homomorphism
$\varphi$  is given not only by an equivalence relation $\Theta$,
defined as usual by $ a \mathrel { \Theta   } b $ if 
$\varphi(a)= \varphi (b)$, but also by a set of further
relations, one additional relation for each relation symbol in  
$\mathscr L$.

\section{Specialization semilattices
``are''  semilattices with a congruence} \label{Cong}

Recall the definition of a specialization semilattice from the 
introduction. It is elementary to see that a specialization semilattice satisfies
\begin{align} 
\label{s7}   
 a \sqsubseteq b  \ \&\  a_1 \sqsubseteq b_1
 & \Rightarrow 
 a \vee a_1 \sqsubseteq b \vee b_1.
\end{align} 
See \cite[Remark 3.5(a)]{mtt}.

If $\mathbf S = (S, \vee)$ is a join semilattice,
we say that a preorder $ \sqsubseteq $  on $S$ is \emph{compatible} if
$(S, \vee, \sqsubseteq )$ is a specialization semilattice, namely,
$ \sqsubseteq $ is coarser than $\leq$ 
(the order induced by $\vee$) and \eqref{s3}
is satisfied.  
As remarked in
\cite[Remark 7.1]{mtt}, the set of all the
 compatible preorders on $\mathbf S$ 
is a complete lattice, with minimum $\leq$,
maximum the total binary relation and meet as intersection
of relations, namely, $a$ and $b$ are related by  
$\bigcap _{i \in I} \sqsubseteq_i$ if 
$a \sqsubseteq_i b $, for every $i \in I$. 

Recall that a \emph{congruence} on a semilattice $\mathbf S$
is an equivalence relation $\Theta$ on $S$ such that  
$ a \mathrel { \Theta} b $ implies
 $ a \vee c \mathrel { \Theta} b \vee c $
for all $a, b,c \in S$. 
This is a special case of a general notion
which applies to arbitrary algebraic structures \cite{Be}.
The set of congruences of some algebraic structure,
in particular, of a semilattice,
is a complete lattice, as well.

\begin{theorem} \label{latt}
Suppose that $\mathbf S = (S, \vee)$ is a join semilattice.

The lattice of compatible preorders on $\mathbf S$ 
is isomorphic, as a complete lattice, 
 to the lattice of congruences on $\mathbf S$.

An isomorphism is given by the function
$\Psi _{_\mathbf S}$ which assigns to a compatible
preorder $ \sqsubseteq $ the congruence
$ \Psi_{_\mathbf S} ( \sqsubseteq ) = \Theta_ \sqsubseteq $  defined by 
$a \mathrel{\Theta_ \sqsubseteq } b$
if both $a \sqsubseteq b$  and $b \sqsubseteq a$. 

The inverse of  
$\Psi_{_\mathbf S}$ is the function
 $\Omega_{_\mathbf S}$  which assigns to a congruence $\Theta$ 
the preorder $ \sqsubseteq _ \Theta $ such that
 $a \sqsubseteq _ \Theta b$ if
$ a \vee b \mathrel { \Theta  } b$.  
 \end{theorem} 

\begin{proof} 
If $ \sqsubseteq $ is a compatible preorder,
then $\Theta_ \sqsubseteq$ is a congruence on $\mathbf S$,
by \eqref{s7} and transitivity and reflexivity of $ \sqsubseteq $.
 Thus the codomain of $\Psi_{_\mathbf S}$
is as wanted. Henceforth we will  drop the subscripts in 
$\Psi_{_\mathbf S}$ and
 $\Omega_{_\mathbf S}$. 

Conversely, let $\Theta$ be a congruence. If $a \leq b$,
that is, $a \vee b = b$,
then $a \vee b \mathrel { \Theta   }  b$ holds, thus
$a \sqsubseteq _ \Theta b$. This shows that
 $\sqsubseteq _ \Theta$ is coarser than $\leq$.
In particular, $\sqsubseteq _ \Theta$ is reflexive.
If $a \sqsubseteq _ \Theta b$ and $b \sqsubseteq _ \Theta c$,
then $a \vee b \mathrel { \Theta   }  b$ and $b \vee c \mathrel { \Theta   }  c$,
thus $a \vee c \mathrel { \Theta   }  a \vee (b \vee c ) =  (a \vee b ) \vee c
\mathrel { \Theta   }  b \vee c \mathrel { \Theta   } c$,
since $\Theta$ is a congruence, hence 
$a \sqsubseteq _ \Theta c$, since $\Theta$ is transitive. This shows that 
$ \sqsubseteq _ \Theta $ is transitive.
If  $a \sqsubseteq _ \Theta b$ and $a_1 \sqsubseteq _ \Theta b$,
then $a \vee b \mathrel { \Theta   }  b$ and $a_1 \vee b \mathrel { \Theta   }  b$,
thus $a \vee a_1 \vee b \mathrel { \Theta   }  b \vee b =b$,
hence $a \vee a_1 \sqsubseteq _ \Theta b$.
We have showed that $\sqsubseteq _ \Theta$ 
is a compatible preorder satisfying \eqref{s3}, hence the codomain of $\Omega$ 
is appropriate, as well.

If $ {\sqsubseteq} = \bigcap _{i \in I} \sqsubseteq_i  $,
then  $a \mathrel{\Theta_ \sqsubseteq } b$
if and only if
 $a \sqsubseteq_i b$  and $b \sqsubseteq_i a$,
for every $i \in I$, if and only if  
$a \mathrel{\Theta_ { \sqsubseteq _i}} b$ for every $i \in I$.
Thus $\Psi (\bigcap _{i \in I} \sqsubseteq_i)
 = \bigcap _{i \in I} \Psi ( \sqsubseteq _i)$,
that is, $\Psi$ is a complete meet homomorphism.  

If we show that $\Omega$ is the inverse of $\Psi$,
then $\Psi$ is bijection, hence an isomorphism.
Note that a complete meet isomorphism
is also a complete join isomorphism.
If $ \sqsubseteq $ is a compatible
preorder and $ \sqsubseteq ^* = \Omega(\Psi( \sqsubseteq ))$,
then $ a \sqsubseteq^* b$ if and only if   
$a \vee b \mathrel { \Theta _ \sqsubseteq   } b$, if and only if  
both $ a \vee b \sqsubseteq b $ and 
 $  b \sqsubseteq a \vee b $. The latter inclusion is always
true, since $ \sqsubseteq $ is coarser than $\leq$;
the former inclusion holds if and only if $a \sqsubseteq b$,
in view of \eqref{s3} and since $ \sqsubseteq $ 
is transitive and coarser then $\leq$. Thus $ a \sqsubseteq^* b$ if and only if 
$ a \sqsubseteq b$, hence $\Psi \circ \Omega$ is the identity.

Conversely, if $\Theta$ is a congruence and 
$\Theta^* =  \Psi(\Omega( \Theta ))$, then
$a \mathrel { \Theta^*} b$ if and only if both
$ a \sqsubseteq _ \Theta b $ and $ b \sqsubseteq _ \Theta a$,
if and only if both
$a \vee b \mathrel { \Theta  } b $  and
$b \vee a \mathrel { \Theta  } a $.
This implies
$ a \mathrel { \Theta  } b$, since $a \vee b = b \vee a$ and $\Theta$ 
is transitive. Conversely, if 
$ a \mathrel { \Theta  } b$, then $a \vee b \mathrel { \Theta  } b $
and $b \vee a \mathrel { \Theta  } a $. 
Thus $a \mathrel { \Theta^*} b$ if and only if 
$a \mathrel { \Theta} b$.
\end{proof}  

\begin{definition} \label{quotie}    
If $\mathbf S$ is a semilattice, a \emph{quotient}
of $\mathbf S$ is a pair $( \varphi, \mathbf T)$, where
$\mathbf T$ is a semilattice and $\varphi$  is a surjective
homomorphism from $\mathbf S$ to $\mathbf T$.
Two quotients $( \varphi, \mathbf T)$ and $( \varphi', \mathbf T')$
of $\mathbf S$ are \emph{isomorphic} if there is
an isomorphism $\chi: \mathbf T \to \mathbf T'$  such that
$  \varphi \circ \chi  = \varphi '$. 
Of course, the above definition applies to any kind of algebraic
or relational structure. 
\end{definition}

From Theorem \ref{latt} and the   Fundamental  Homomorphism  Theorem
\cite[Theorem 1.22]{Be} we immediately get the following 
corollary.

\begin{corollary} \label{quot}
Let $\mathbf S$ be a  semilattice and consider 
the function $\Gamma$ which sends a compatible preorder
$ \sqsubseteq $ to the pair 
$( \pi_{ \sqsubseteq }, \mathbf S/\Theta _{ \sqsubseteq }  )$,
where $\Theta _{ \sqsubseteq }$ is defined in Theorem \ref{latt} and
$\pi_{ \sqsubseteq }$ is the canonical projection
 from $S$ to $S/\Theta _{ \sqsubseteq }$.    

Then $\Gamma$ is a bijection from the set of  compatible preorders
on $\mathbf S$ to the set of equivalence classes of quotients of $\mathbf S$ 
up to isomorphism. 
 \end{corollary}   

As another corollary of Theorem \ref{latt}, we get
the alternative representation theorem for 
specialization semilattices mentioned in the introduction. 

\begin{corollary} \label{rep}
Suppose that $(S, \vee_S)$
is a  semilattice and $ \sqsubseteq $ is a binary
relation on $S$.

Then $(S, \vee_S, \sqsubseteq )$
is a specialization semilattice if and only if
 there are a semilattice $(T, \vee_T)$ 
and a semilattice homomorphism
$\varphi: (S, \vee_S) \to (T, \vee_T)$ such that 
$a \sqsubseteq b$ if and only if $\varphi(a) \leq_T \varphi (b) $,
for all $a,b \in S$.  
Namely,   $ \sqsubseteq $ is $ \sqsubseteq _ \varphi  $,
as given by \eqref{sqf}.  
 \end{corollary} 

\begin{proof}   
The ``if'' condition is straightforward. 
Conversely, if $(S, \vee_S, \sqsubseteq )$
is a specialization semilattice,
then $\Theta_ \sqsubseteq $ is a congruence of
$(S, \vee_S)$, by Theorem \ref{latt}.
If $(T, \vee_T)$ is the quotient of 
$(S, \vee_S)$ modulo $\Theta_ \sqsubseteq $
with
$\varphi$  is the canonical projection,
then, for all $a,b \in S$, we have $\varphi(a) \leq_T \varphi (b) $
if and only if 
$\varphi(a \vee_S b) = \varphi(a) \vee_T \varphi(b) = \varphi(b)  $,  
if and only if 
  $a \vee_S b \mathrel { \Theta _ \sqsubseteq } b$.
Thus  $ {\sqsubseteq _ \varphi}
  = \Omega(\Psi( \sqsubseteq )) ={ \sqsubseteq }$
by Theorem \ref{latt}.
\end{proof}

We now reformulate the above results in a categorical
framework.
We refer to \cite{AHS} for basic notions about
categories.
 
Specialization semilattices form a category
with the standard model theoretical notion of homomorphism,
namely, a homomorphism between two specialization semilattices
is a semilattice homomorphism $\varphi$  such that $ a \sqsubseteq b $
implies $ \varphi (a) \sqsubseteq \varphi (b)$.
The category $\mathcal S \mathfrak{Con}$ of \emph{semilattices with a congruence}
 has objects those triplets $(S, \vee, \Theta )$    
such that $(S, \vee)$ is a semilattice and $\Theta$ 
is a congruence on $(S, \vee)$. A morphisms $\varphi$  in
$\mathcal S \mathfrak{Con}$
from $(S, \vee_S, \Theta )$ to $(T, \vee_T, \Gamma  )$
is a  semilattice homomorphism such that $a \mathrel { \Theta  } b $ 
implies $ \varphi (a) \mathrel { \Gamma   } \varphi (b) $.
Of course, semilattices play no special role in the above definition:
for every variety $\mathcal V$ of algebraic structures,
exactly in the same way,
we can define the category $\mathcal V \mathfrak{Con}$
 of algebras in $\mathcal V$  with a congruence.

\begin{corollary} \label{corcat}
The category of specialization semilattices
is isomorphic to the category
of semilattices with a congruence.

The correspondence between objects is given by 
the functions $\Psi _{_\mathbf S}$ from Theorem \ref{latt}
and  by the identity for morphisms.  
 \end{corollary} 

\begin{proof}
In view of Theorem \ref{latt}, it is enough to show that
some function  $\varphi$  is a morphism for specialization semilattices
if and only if  $\varphi$  is a morphism for 
the corresponding semilattices
with a congruence.

Indeed, if $\varphi$  is a morphism
from $(S, \vee_S, \sqsubseteq _S)$ to $(T, \vee_T, \sqsubseteq _T  )$
and $a  \mathrel { \Theta _{\sqsubseteq _S} } b $,
then $a \sqsubseteq_S b $ and $ b \sqsubseteq_S a $,
thus  $ \varphi (a) \sqsubseteq_T \varphi ( b) $
 and $ \varphi (b) \sqsubseteq_T \varphi (a) $,
hence $ \varphi (a)  \mathrel { \Theta _{\sqsubseteq _T} } \varphi (b) $.

Conversely, if $\varphi$  is a morphism
from $(S, \vee_S, \Theta)$ to $(T, \vee_T, \Gamma  )$
and $ a \sqsubseteq _ \Theta b  $, then
$a \vee b \mathrel { \Theta} b $, thus 
$ \varphi (a) \vee \varphi ( b) =
 \varphi (a \vee b) \mathrel {\Gamma} \varphi (b) $,
since $\varphi$  is a morphism
of semilattices with a congruence, in particular, a semilattice homomorphism.   
Hence $ \varphi (a) \sqsubseteq _ \Gamma  \varphi (b)  $.
 \end{proof}  

By a \emph{concrete} category we mean a \emph{construct}
in the terminology from \cite{AHS}, namely, a concrete category
over \textbf{Set}.    
If $\mathfrak C$ is a (concrete) category,
the  category $\mathfrak {CEpi}$ of  $\mathfrak C$-epimorphisms
(the  category $\mathfrak {CSurj}$ of  $\mathfrak C$-surjections)
is the category whose objects are epimorphisms of 
$\mathfrak C$ (surjective morphisms of 
$\mathfrak C$). The morphisms of $\mathfrak {CEpi}$
($\mathfrak {CSurj}$)
are commuting squares
of morphisms of $\mathfrak C$. 
Thus $\mathfrak {CEpi}$ and $\mathfrak {CSurj}$ are  full  subcategories
of the arrow category of $\mathfrak C$;
 in $\mathfrak {CEpi}$  ($\mathfrak {CSurj}$) only 
epimorphisms (surjective morphisms) are considered as ``arrows'', but 
the other morphisms in the commuting squares are not
necessarily assumed to be epi or surjective.

\begin{remark} \label{cattt}    
Congruences in general algebraic systems
 are just an internal description
of homomorphic images.
Restated from a categorical point of view, if $\mathcal V$ is a variety, 
then the category of algebras in $\mathcal V$  with a congruence
is equivalent to the  
 category $\mathcal V \mathfrak {Surj}$ of 
 $\mathcal V$-surjections.
Indeed, consider the functor $F$ which associates to an
algebra with a congruence $(\mathbf A, \Theta )$ 
the canonical quotient morphism
$\pi_\Theta : \mathbf A \to \mathbf A/ \Theta $.
If $\varphi:A \to B$ is a function,
then   $\varphi $ is a morphism
of $\mathcal V\mathfrak{Con}$ from 
$  (\mathbf A, \Theta ) $ to $  (\mathbf B, \Gamma  ) $
if and only if $\varphi$   is a morphism in $\mathcal V$ 
 and there exists a $\mathcal V$-morphism $\hat \varphi $ such that  
 $  \varphi \circ \pi_\Gamma = \pi_\Theta \circ \hat \varphi  $, thus
necessarily $\hat \varphi $ is defined by 
$ a /\Theta   \mapsto \varphi(a) / \Gamma$. 
In fact, this assignment provides a good definition exactly when
$\varphi$  is a $\mathcal V\mathfrak{Con}$-morphism from 
$  (\mathbf A, \Theta ) $ to $  (\mathbf B, \Gamma  ) $.
Thus if we set  $F( \varphi ) = ( \varphi , \hat \varphi ) $  on morphisms, then 
$F$ is full and faithful from $\mathcal V\mathfrak{Con}$ to 
$\mathcal V \mathfrak {Surj}$.  If 
$ \psi : \mathbf A \to \mathbf  B$ is surjective,
then $\psi$ is isomorphic to the canonical
projection
$\pi_\Theta : \mathbf A \to \mathbf A/ \Theta $
in  $\mathcal V\mathfrak {Surj}$, where $\Theta$ 
  is the kernel of $\psi$, namely
$\Theta =\{ \, (a,b) \mid  \psi(a) =\psi(b) \,\}$. 
This means that $F$ is \emph{isomorphism-dense}
\cite[Definition 3.33]{AHS}.  
 Thus $F$ witnesses that 
the category of algebras in $\mathcal V$  with a congruence
is equivalent to the 
category of  $\mathcal V$-surjections.

As we will see in Section \ref{gen},
the situation is more delicate in the context of
relational structures. For arbitrary categories, additional problems might
arise;
see, e.~g., the section ``Quotient objects'' in  \cite{AHS}.
\end{remark}

The above considerations apply, in particular, to 
the variety of semilattices. Noticing that in the category
of semilattices epimorphisms are exactly surjective morphisms \cite{KMP},
from Corollary \ref{corcat}  we get the following corollary.

\begin{corollary} \label{equiv}
The category of specialization semilattices is equivalent
to the   category of semilattice epimorphisms.
 \end{corollary} 

\begin{remark} \label{diff1}   
There is a subtle difference
between Corollary \ref{equiv} and
Corollary \ref{quot}. 
For example, let $\mathbf S$ be a semilattice
with five elements $a_1, a_2, b_1, b_2, c$,
with  $a_1 < a_2 < c$,
$ b_1 < b_2 <c$ and 
$a_i \vee b_j = c$, for $i, j \in \{ 1,2\}$.
Let $\Theta$ (resp. $\Theta'$)  be the equivalence relation 
such that  $a_1$ and $a_2$ are related
(resp., $b_1$ and $b_2$ are related) and no
other pair of distinct elements are related.
Then the corresponding quotients $( \pi, \mathbf T)$ 
and $(\pi', \mathbf T')$ are not isomorphic. 
Indeed, the only semilattice isomorphism $\chi$ from
$\mathbf T$ to $\mathbf T'$ sends the class 
$\{ a_1, a_2\}$ to the class $\{ b_1, b_2\}$,
but then   $\pi \circ \chi = \pi'$ fails.

On the other hand, the objects 
$ \mathbf S \stackrel{\pi}{\to}  \mathbf T $   
and 
$ \mathbf S \stackrel{\pi'}{\to} \mathbf T' $
are isomorphic in $\mathcal S\mathfrak{Epi}$.   
 \end{remark}

\section{An internal characterization of surjective
homomorphisms between relational systems} \label{gen}

We now generalize some arguments from the previous section.
We begin with an example.

\begin{example} \label{sppo}   
A \emph{specialization poset} is a partially ordered set together
with a preorder coarser than the order.
If $X$ is a topological space,
then  $(\mathcal P(X), \subseteq , \sqsubseteq )$ 
is a specialization poset.
As recalled in the introduction, 
and parallel to the case of specialization semilattices, 
every specialization poset can be embedded into
some ``topological'' specialization poset as above;
moreover, continuous functions correspond
exactly to specialization homomorphisms 
\cite[Propositions 5.10  and 2.4]{mtt}. 

If $\mathbf P= (P, \leq_P)$ and 
$\mathbf Q =(Q, \leq_Q)$ are posets, 
$\varphi : \mathbf P \to \mathbf Q$    
is an order preserving function and,
for $a,b \in P$,  we set 
$a \sqsubseteq_ \varphi  b$ if $\varphi(a) \leq_Q \varphi (b)$,
as in \eqref{sqf}, 
then   $(P, \leq_P, \sqsubseteq _ \varphi )$
is a specialization poset. 
As in the previous section, every specialization poset can be
constructed in this way, but in this case the argument is much easier.
Indeed, if $ \sqsubseteq $ is a preorder on some set $P$,
then, classically, $ \sqsubseteq $ induces a partial order on $P/ {\sim}$,
where $\sim$ is the equivalence relation on $P$ defined
by $a \sim b$ if both $ a \sqsubseteq b$ and $b \sqsubseteq a$.    
If $\leq$ is an order on $P$ and $\leq$ is finer than $ \sqsubseteq $,
then the projection $\pi$ which sends $a$ to $ a/ {\sim}$ is order
preserving from $(P, \leq)$ to  $(P/ {\sim}, \sqsubseteq / {\sim})$.
Then $\sqsubseteq_ \pi$ coincides with the original $ \sqsubseteq $,
hence we have got the desired representation.
\end{example} 

By using some arguments  from Remark \ref{cattt},
we then get the following proposition.  

\begin{proposition} \label{equivb}
The category of specialization posets is equivalent
to the   category $\mathcal{POS}\mathfrak{Surj}$, whose objects are  
 surjective order preserving
maps between posets. Morphisms of 
$\mathcal{POS}\mathfrak{Surj}$
are commuting squares of order preserving maps.
 \end{proposition} 

\begin{remark} \label{diff}    
Notice a difference in comparison with the previous section.
Here we cannot do simply with an equivalence relation.
Consider the poset $\mathbf P$ with just two incomparable elements
$a$ and  $b$. Up to isomorphism,
we have four surjective order preserving functions from $\mathbf P$
to some other poset: the 
isomorphism, the quotient which identifies $a$ and $b$
and the two bijections onto a poset in which $a$ and  $b$  
become comparable\footnote{These last two are better considered 
to be distinct,
since we  consider surjective functions up to 
isomorphism, but keeping $\mathbf P$ fixed, thus 
the elements of $P$ are considered to ``have a name''. However, the
remark in the present footnote is irrelevant for   what follows.}.
On the other hand, there are only two equivalence relations
on a set with two elements, hence surjective order preserving maps
cannot be characterized just by equivalence relations.
See \cite{W} for a survey of what can be done with just an
equivalence relation. Note also that the image of a poset under
some surjective function need not be transitive, neither antisymmetric.
 See again  \cite{W}.
\end{remark}

We now see that Example \ref{sppo}  can be generalized 
to arbitrary models.
For the rest of this section we assume that the reader has
some familiarity with the basic notions of model theory.
See e.~g. \cite{H}.

\begin{definition} \label{def}
(a) Suppose that $\mathbf A$ and $\mathbf  B$  are models
for the same language $\mathscr L$ 
and $\varphi: \mathbf A \to \mathbf  B$ is a surjective homomorphism.
We say that $\varphi$  is \emph{into classes}
if there is some equivalence relation $\Theta$ on $A$ 
such that $B= \{ \, a/ \Theta  \mid a \in A \, \} $, as a set,
and moreover $\varphi(a)= a/ \Theta$, for every $a \in A$.
Of course, if this is the case, then $\Theta = \{ \, (a, a_1) \mid
\varphi (a) = \varphi (a_1)  \, \} $.     
Clearly, every  surjective homomorphism 
is isomorphic (in the sense of Definition \ref{quotie}) to some  surjective homomorphism into classes. 

(b)  Suppose that $\mathbf A$ is a model
for the  language $\mathscr L$.
An expansion $\mathbf A^{\myplus}$  of $\mathbf A$ is
\emph{appropriate for describing a surjective homomorphic
image}, for short, simply \emph{appropriate}, if
the following conditions hold.

(1) $\mathbf A^{\myplus}$ is a model for the language
$\mathscr L^{\myplus} $, where $\mathscr L^{\myplus} $ 
is obtained from
$\mathscr L$ by adding,
for every relation symbol $R$ of $\mathscr L$,
 a new relation symbol
$R^{\myplus}$ of the same arity.
This is  applied also to the
equality symbol $\equiv$,
and we will write $\Theta$ for $\equiv^{\myplus}$.
Notice that $\mathscr L$ is allowed to 
contain constant and function symbols, but no
 new constant nor function symbol is  added to
$\mathscr L^{\myplus} $.
In particular, if $\mathscr L$ contains no relation symbol,
apart from equality, 
then $\Theta$ is the only additional relation symbol in
$\mathscr L^{\myplus} $.

(2)
The interpretation of
$\Theta$ in $\mathbf A^{\myplus}$ 
 is an equivalence relation and, furthermore, the 
universal closures of the following
formulas are satisfied in $\mathbf A^{\myplus}$
\begin{align}
\label{axx0}
 &R (x_1, \dots, x_n) \Rightarrow 
 R^{\myplus}(x_1, \dots, x_n)
\\ 
\label{axx1}  
 & x_1 \mathrel { \Theta  } y_1 \wedge \dots \wedge
 x_n \mathrel { \Theta  } y_n \Rightarrow 
 f (x_1, \dots, x_n) \mathrel { \Theta  } f (y_1, \dots, y_n),
\\
\label{axx2}  
& x_1 \mathrel { \Theta  } y_1 \wedge \dots \wedge
 x_n \mathrel { \Theta  } y_n  \wedge 
 R^{\myplus} (x_1, \dots, x_n) \Rightarrow  R^{\myplus} (y_1, \dots, y_n),
 \end{align}
for every function symbol $f$
and every relation symbol $R$ in $\mathscr L$,
where we implicitly assume that indexes are appropriate
for arities.
 \end{definition}   

In the next proposition the symbols for the 
$R^{\myplus}$s are kept fixed; otherwise, the result holds
just up to renamings. 

\begin{proposition} \label{chara}
Let $\mathbf A$ be a model
for the  language $\mathscr L$.
There is a bijection between the set of 
appropriate
expansions $\mathbf A^{\myplus}$  of $\mathbf A$
and the set of  surjective homomorphisms into classes
with domain $\mathbf A$ or, put in an equivalent way,
 the set of equivalence classes of
surjective homomorphisms
with domain $\mathbf A$, up to isomorphism
as in Definition  \ref{quotie}.
\end{proposition} 

 \begin{proof}
If $\varphi: \mathbf A \to \mathbf  B$ is a surjective homomorphism,
not necessarily into classes, expand $\mathbf A$ 
to an $\mathscr L^{\myplus}$ model 
by interpreting $\Theta$ as  $\{ \, (a, a_1) \mid
\varphi (a) = \varphi (a_1)  \, \} $ and, for each relation symbol $R$ in
$\mathscr L$,  letting
$ R^{\myplus} (a_1, \dots, a_n)$ hold in $\mathbf A^{\myplus}$
if  $ R( \varphi (a_1), \dots, \varphi ( a_n))$ holds in $\mathbf  B$.
Then $\mathbf A^{\myplus}$  is appropriate; indeed, 
\eqref{axx0} and \eqref{axx1} follow from the assumption that
$\varphi$  is a homomorphism, while \eqref{axx2} holds by construction.

Conversely,   if $\mathbf A^{\myplus}$  is an appropriate
expansion of $\mathbf A$, let $B= B / \Theta$
and, for every relation symbol $\mathcal R \in \mathscr L$,
 let   $ R ( a_1 / \Theta,
 \dots, a_n/ \Theta)$ hold in $\mathbf  B$
if $ R^{\myplus} (a_1, \dots, \allowbreak  a_n)$ holds in $\mathbf A^{\myplus}$.
This is a good definition in view of \eqref{axx2}.
As standard, for function symbols,
 let $ f ( a_1 / \Theta,
 \dots, a_n/ \Theta) = f ( a_1 ,
 \dots, a_n) / \Theta $ in $\mathbf  B$.
This is well-defined because of \eqref{axx1}. 
The projection $\pi$ sending $a \in A$ to
$a/\Theta$ is then a surjective homomorphisms into classes,
by construction and \eqref{axx0}.

 If $\varphi$  in the first paragraph of the proof is actually into classes,
then the above constructions are one the inverse of the other,
hence they provide a bijection as desired.
The last statement is immediate from
the last comment in Definition \ref{def}(a).
\end{proof} 

If $\mathscr L$ contains no relation symbol apart from equality,
then the only additional relation in $\mathbf A^{\myplus}$ is $\Theta$.
In this case $\Theta$ is interpreted as a congruence on the 
algebraic structure $\mathbf A$, hence Definition \ref{def}(b)
generalizes the notion of a congruence and Proposition \ref{chara}  
generalizes the basic property of a congruence.
Note also that, in the special case of posets,
 the content of Proposition \ref{chara} does not correspond
exactly to the content of Example \ref{sppo}, since, in this special case
of Proposition \ref{chara}, $\leq^{\myplus}$ is not assumed
to be transitive, namely, $\mathbf  B$ does not necessarily become
a poset. Compare the last sentence in Remark \ref{diff}. 

Proposition \ref{chara} is essentially a reformulation of
known ideas; see  \cite[Section 1.4]{G}, whose results and
methods, in our modest
opinion, have not yet received the full attention they deserve.
Indeed, very special cases have been subsequently and independently
discovered.
 The special case of graphs has been  dealt with in \cite{BHV},
with many applications. The special case
of ordered algebras in which operations are either monotone
or antitone at each place has been dealt with in \cite{P}
(\cite{P} deals only with $\leq^{\myplus}$, in the present notation,
but, under the assumptions, $\Theta$ can be defined in terms of  $\leq^{\myplus}$).  
 Note that, in the case of posets, 
if one wants to deal with only an equivalence relation
 and no
auxiliary relation, the situation appears much more involved and with
 many facets \cite{W}. 
Many more  references can be found in the quoted works.

\end{document}